\documentclass{amsart}
\usepackage{graphicx} 

\usepackage{amsmath}
\usepackage{amssymb}
\usepackage{amsthm}
\usepackage[normalem]{ulem}
\usepackage{dsfont}

\usepackage{esint}

\usepackage{mathrsfs}

\theoremstyle{theorem}
\newtheorem{theorem}{Theorem}[section]

\theoremstyle{definition}
\newtheorem{definition}[theorem]{Definition}

\theoremstyle{definition}

\theoremstyle{definition}

\theoremstyle{theorem}
\newtheorem{lemma}[theorem]{Lemma}

\theoremstyle{theorem}
\newtheorem{proposition}[theorem]{Proposition}

\theoremstyle{theorem}

\theoremstyle{corollary}
\newtheorem{corollary}[theorem]{Corollary}

\theoremstyle{remark}
\newtheorem{remark}[theorem]{Remark}

\DeclareMathOperator{\tr}{tr}
\DeclareMathOperator{\rank}{rank}

\usepackage{polynom}

\DeclareMathOperator{\End}{End}

\DeclareMathOperator{\Hom}{Hom}

\DeclareMathOperator{\ind}{ind}

\DeclareMathOperator{\Pin}{Pin}
\DeclareMathOperator{\Spin}{Spin}

\DeclareMathOperator{\delbar}{\overline{\partial}}

\usepackage{comment}
\usepackage{appendix}

\usepackage{tikz-cd}
\usetikzlibrary{matrix,arrows,decorations.pathmorphing}

\usepackage{hyperref}

\title[Vector bundles on spin K\"ahler surfaces]{Vector bundles on spin K\"ahler surfaces and the multi-spinor Seiberg-Witten equations}
\author{Minh Lam Nguyen}
\address{\parbox{\linewidth}{Department of Mathematics, The University of Texas at Tyler, Texas, 75799}}
\email{minhnguyen@uttyler.edu}
\author{Ollie Thakar}
\address{\parbox{\linewidth}{Department of Mathematics, Harvard University, Massachusetts, 02138}}
\email{othakar@math.harvard.edu}

\begin{document}

\maketitle
\begin{abstract}
    We compute new lower bounds on the second Chern class of a vector bundle $\mathscr{E}$ on a spin K\"ahler surface $X$ that satisfies $c_1(\mathscr{E}) = 0$ and $H^2(\End_0\mathscr{E}) = 0.$ When $X$ is a regular elliptic surface and $\mathscr{E}$ has rank 2, our bound is sharp. Our technique is to use Furuta's finite-dimensional approximation on the multi-spinor Seiberg-Witten equations, together with their algebro-geometric interpretation by the second named author. We also use these techniques to prove an existence theorem for singular Fueter sections and multi-spinor Seiberg-Witten solutions that holds for every Riemannian metric on a large class of 4-manifolds.
\end{abstract}

\section{Introduction}

In this paper, we use the multi-spinor Seiberg-Witten equations and Furuta's technique of finite-dimensional approximation to study stable vector bundles on spin K\"ahler surfaces. Given a K\"ahler surface $X$, we ask, if there exists a vector bundle $\mathscr{E}\to X$ with $c_1(\mathscr{E}) = 0$ and $H^2(\End_0\mathscr{E}) = 0,$ what values of $c_2(\mathscr{E})$ are possible? If $\mathscr{E}$ is stable in the sense of Mumford-Takemoto, the condition that $H^2(\End_0\mathscr{E}) = 0$ guarantees that in a neighborhood of $\mathscr{E}$ the moduli space of stable bundles on $X$ is smooth and its dimension equals its virtual dimension. This existence question and its variations have spawned much work in gauge theory and algebraic geometry \cite{Taubes, Gieseker, Friedman-elliptic, donaldsonpoly, OGrady}, largely catalyzed by Donaldson's work on the Yang-Mills equations on algebraic surfaces. To state our results, we must introduce two definitions:

\begin{definition}[Friedman {\cite[Chapter 6, Definition 16]{Friedman}}]\label{good}
    For a vector bundle $\mathscr{E}$, let $\End_0\mathscr{E}$ denote the bundle of traceless endomorphisms. A holomorphic vector bundle $\mathscr{E}\to X$ over a complex manifold $X$ is \emph{good} if $H^2(\End_0\mathscr{E}) = 0.$
\end{definition}

\begin{definition}
    Let $(X, \omega)$ be a compact K\"ahler surface. For each positive integer $r,$ define $n_r(X)$ (resp. $n_r(X, \omega)$) to be the smallest value of $n\in\mathbb{Z}$ for which there exists a good (resp. good and $\omega$-stable) holomorphic vector bundle $\mathscr{E}\to X$ with: $$\text{rank}~ \mathscr{E} = r, ~c_1(\mathscr{E}) = 0, ~c_2(\mathscr{E}) = n.$$
\end{definition}

By the work of Taubes, $n_2(X)$ is finite \cite[Theorem 1.1]{Taubes}; this was extended to $n_r(X)$ for all $r\geq 2$ by O'Grady \cite[Corollary B']{OGrady}. Many authors have found upper bounds on these invariants \cite{Gieseker, Friedman-elliptic}. We will prove the following \emph{lower bound} on $n_2(X)$, where $\sigma(X)$ is the signature of $X$ (see Corollary \ref{su(2)}):

\begin{theorem}
    Let $X$ be a compact spin K\"ahler surface. We have: $$n_2(X) \geq -\frac{\sigma(X)}{4} = 2\chi(\mathscr{O}_X) - \frac14c_1^2(X).$$
\end{theorem}

This bound is sharp for a large class of surfaces, allowing us to compute $n_2(X)$ explicitly for this class (as far as we know, there are very few varieties $X$ for which $n_2(X)$ has been computed explicitly.) We follow the conventions of algebraic surfaces and use $p_g(X)$ to denote the geometric genus $h^2(\mathscr{O}_X).$

\begin{theorem}[will be proved as Theorem \ref{elliptic}]
    If $(X,H)$ is a polarized elliptic surface with $h^1(\mathscr{O}_X) = 0$ and $K_X$ even, where $H$ is an ample line bundle over $X$, then $$n_2(X) = n_2(X, H) = 2p_g(X) + 2.$$
\end{theorem}

\begin{remark}
    The condition that $X$ is spin in both of these theorems is necessary, and cannot even be weakened to the condition that $X$ is minimal. For, consider the Dolgachev surface $Z$ constructed by applying logarithmic transforms to two distinct smooth fibers of the elliptic fibration $\mathbb{CP}^2\#9\overline{\mathbb{CP}}^2\to \mathbb{CP}^1$, so that the resulting multiple fibers have multiplicities 2 and 3, respectively. This is the surface considered by Donaldson in his disproof of the 4-dimensional $h$-cobordism conjecture \cite{Donaldsonelliptic}. The surface $Z$ is a minimal projective elliptic surface with $\pi_1(Z) = 1$, so in particular $h^1(\mathscr{O}_Z) = 0$. By \cite[Theorem 3.21]{Donaldsonelliptic}, there exists a polarization $H$ and an $H$-stable good rank 2 holomorphic vector bundle $\mathscr{E}\to Z$ with $c_1(\mathscr{E}) = 0,~ c_2(\mathscr{E}) = 1.$ But, $\sigma(Z) = -8,$ so $Z$ is not spin by Rokhlin's Theorem, and we see: $$n_2(Z, H)\leq 1<2= \frac{-\sigma(Z)}{4}.$$ 
\end{remark}

\begin{remark}
    These new lower bounds on $n_2(X)$ allow us to conclude in multiple examples that Donaldson polynomial invariants of an algebraic surface can be nonzero even when the corresponding anti-self-dual Yang-Mills moduli spaces consist of only obstructed connections for all K\"ahler metrics; see Remark \ref{donaldsonpoly}.
\end{remark}

\begin{remark}
    A compact K\"ahler manifold $X$ is spin if and only if its canonical class $K_X$ is even; we will therefore interchange freely between these two notions.
\end{remark}

We also prove lower bounds for the second Chern classes of stable good bundles $\mathscr{E}\to X$ of arbitrary rank that admit real or quaternionic structures (see Corollary \ref{realstable} and Corollary \ref{sp}.) The most general result is as follows:

\begin{theorem}[will be proved as Theorem \ref{selfdual}]
    If $X$ is a spin K\"ahler surface and $\mathscr{E}\to X$ is a good holomorphic vector bundle of rank $r$ which is holomorphically isomorphic to its dual bundle, then: $$c_2(\mathscr{E}) \geq -r\frac{\sigma(X)}{8} = r\left(\chi(\mathscr{O}_X) - \frac18 c_1^2(X)\right)$$
\end{theorem}

\subsection{Multi-Spinor Seiberg-Witten Theory}

The technique we use combines two ideas which are novel to the study of algebraic surfaces. So, we summarize the way we incorporate these ideas below. 

The first is the study of the multi-spinor Seiberg-Witten equations on algebraic surfaces; this was initiated by Bryan and Wentworth \cite{Bryan-Wentworth} and generalized in the setting we use here, where the data of the equations includes an auxiliary vector bundle and connection, by the second named author \cite{Thakar26}. The key idea is that the moduli space of solutions of the multi-spinor Seiberg-Witten equations modulo gauge is compact whenever the auxiliary connection is an anti-self-dual Yang-Mills connection whose deformation complex has vanishing second cohomology (see \cite[Theorem 4.3]{Thakar26} and Theorem \ref{compactness} below.) Unlike their famous single-spinor counterparts, these equations have not so far been successful at proving theorems in differential or algebraic geometry (with the exception of the obscure \cite[Theorem 6.10]{Thakar26}.) So, this paper represents one of the first applications of these equations to prove a purely geometric theorem.

The second idea is to imitate Furuta's proof of the $\frac{10}8+2$ Theorem in 4-manifold topology, which we sketch:

\begin{theorem}[{\cite[Theorem 1]{Furuta}}]\label{furuta}
    If $X$ is a smooth closed spin 4-manifold with indefinite intersection form, then $$b_2(X) \geq \frac54|\sigma(X)| + 2.$$
\end{theorem}

His proof proceeds by considering the Seiberg-Witten equations as a map of Hilbert spaces $\mathcal{F}:V\to W$. When $X$ is spin, there exist actions of the group $\Pin(2)$ on $V$ and $W$ for which $\mathcal{F}$ is equivariant. Using the compactness of the space of solutions, he shows that there exist finite-dimensional $\Pin(2)$-submodules $V_\lambda\subset V$ and $W_\lambda\subset W$ and constructs out of $\mathcal{F}$ a $\Pin(2)$-equivariant map of $\Pin(2)$ unit spheres $S_1(V_\lambda)\to S_1(W_\lambda)$. Using the equivariant Atiyah-Singer Index Theorem, he computes the decomposition of those modules into irreducible representations satisfying certain relations. The multiplicities of those relations are given by topological invariants of $X$. Using equivariant $K$-theory, he finds an inequality that the multiplicities of those irreducible representations must satisfy. That inequality, interpreted as an inequality of those topological invariants, is precisely the desired bound.

Our observation, which originates from the first named author's work on the Rarita-Schwinger-Seiberg-Witten equations \cite{Nguyen}, is that Furuta's use of finite-dimensional approximation and equivariant $K$-theory may be applied to equations other than the Seiberg-Witten equations provided they have compact moduli spaces of solutions and are $\Pin(2)$-equivariant. The multi-spinor Seiberg-Witten equations are $\Pin(2)$-equivariant whenever the base manifold is spin and the auxiliary connection is real. In such a setting, Furuta's machinery then tells us that whenever the moduli space of solutions is compact, we get an inequality among the topological quantities of $X$ and the auxiliary vector bundle (Theorem \ref{Furuta}). Combining this with our algebro-geometric condition for compactness carefully, we arrive at our desired bound on $n_2(X)$ (Corollary \ref{su(2)}) and other useful inequalities.

\subsection{Applications to geometric analysis}

Finite-dimensional approximation for multi-spinor Seiberg-Witten equations also leads to a result in geometric analysis of independent interest, so although it is logically separate from our algebro-geometric theorems we include it in this paper for record. Following the strategy of the first named author in \cite{Nguyen}, we show an existence theorem for multi-spinor Seiberg-Witten solutions and their rescaled limits (see Theorem \ref{existence}.) The idea is that when the multi-spinor Seiberg-Witten moduli space is compact for a fixed choice of Riemannian metric and auxiliary connection, Theorem \ref{Furuta} provides us an inequality relating topological quantities of our manifold and auxiliary bundle. So, when that inequality fails, which is a purely topological condition, the multi-spinor Seiberg-Witten moduli space must be non-compact (in particular, non-empty) for \emph{every} Riemannian metric and auxiliary connection. In addition, Taubes has shown that a non-convergent sequence of multi-spinor Seiberg-Witten solutions, upon rescaling, converges to a \emph{singular Fueter section}, an interesting geometric-analytic object in its own right \cite{Taubes16}. We will therefore also be able to show the existence of such singular Fueter sections under the same topological assumptions.

\subsection{Outline of Paper}

Sections 2 and 3 will be devoted to applying Furuta's machinery to the multi-spinor Seiberg-Witten equations. In Section 2, we will review the equations and verify their $\Pin(2)$-symmetry. In Section 3, we will construct the finite-dimensional approximation and use it to extract a gauge-theoretic inequality (Theorem \ref{Furuta}) analogous to that in Theorem \ref{furuta}.

Section 4 is devoted to applications to algebraic geometry. The first four subsections will comprise the compactness theorem we use to translate our gauge-theoretic inequality into a purely algebro-geometric result, as well as indirect applications of this inequality to conclude the main result of our paper, Corollary \ref{su(2)}. We also get bounds for vector bundles with real and quaternionic structures. The remaining sections compare our bound Corollary \ref{su(2)} to existing bounds and compute it in several examples, notably computing $n_2(X)$ exactly for a large class of elliptic surfaces (Section 4.6 and Theorem \ref{elliptic}) and improving on known lower bounds for $n_2(X)$ for branched covers of $\mathbb{CP}^2$ (Section 4.7).

Finally, Section 5 discusses the applications of Theorem \ref{Furuta} to geometric analysis, namely an existence theorem (Theorem \ref{existence}) for multi-spinor Seiberg-Witten solutions and singular Fueter sections. It is independent of Section 4.

\subsection{AI Disclosure} No LLMs were used in the process of producing this paper.

\subsection{Acknowledgments} MN would like to thank Greg Parker for helpful discussions, and both authors would like to thank him for detailed feedback on an earlier draft. OT would like to thank Dylan Galt for removing an assumption in Theorem \ref{compactness}, and his advisor Peter Kronheimer for sharing invaluable wisdom and support. OT is supported by the Simons Collaboration on New Structures in Low-Dimensional Topology.

\section{Multi-Monopole Map}

In this section, we will define the moduli space of solutions to the multi-spinor Seiberg-Witten equations and prove basic properties thereof.

\subsection{Geometric Set-Up}

Following \cite{Thakar26}, we define the multi-spinor Seiberg-Witten map. The multi-spinor Seiberg-Witten map are gauge-theoretic equations on a closed, smooth, oriented 4-manifold $X$ with the following additional data:

\begin{itemize}
    \item A Riemannian metric $g$
    \item A $SU(r)$-bundle $E\to X$ with $c_2(E)\geq0.$
    \item A $SU(r)$-compatible connection $B\in\mathscr{A}(E)$
    \item A spin$^c$-structure $\mathfrak{s}$ with spinor bundles $W^\pm\to X,$ determinant bundle $S,$ and with Clifford multiplication $\rho:T^*X\otimes \mathbb{C}\to\Hom(W^\pm, W^\mp)$
\end{itemize}

The Seiberg-Witten map will take in a unitary $A$ on the determinant bundle $S$ and a section $\Phi\in\Omega^0(W^+\otimes E)$ which will be referred to as a \emph{twisted spinor.} There is a unique spin$^c$ connection on $W^\pm$ which induces the connection $A$ on the determinant bundle, and we will also label this connection as $A.$ The Seiberg-Witten map $\mathcal{F}:\Omega^0(W^+\otimes E)\times \mathcal{A}(S)\to \Gamma_X(W^-\otimes E)\times \Omega^+(i\mathbb{R})$ is: \begin{equation}\label{msw}\mathcal{S}(A, \Phi) = \begin{pmatrix}
    D_{A\otimes B}^+\Phi \\
    F_{A}^+ -\rho^{-1}\mu(\Phi, \Phi), \\
\end{pmatrix}\end{equation} where $\sigma:(W^+\otimes E)\times\overline{(W^+\otimes E)}\to \mathfrak{sl}(W^+)$ is a bilinear map given by: $$\mu(\Phi, \Psi) := \tr_E(1-\tr_{W^+})(\Phi\otimes \Psi^*).$$

Solutions to these equations are preserved by the action of the gauge group $\mathcal{G} = \text{Map}(X, S^1)$ given as $u\cdot (A, \Phi) = (A - u^{-1}du, u\cdot \Phi).$ The moduli space of such solutions, modulo unitary gauge, is to be called $\mathcal{M}(\mathfrak{s}, g, B).$ 

\subsection{$\Pin(2)$ Symmetry}

\begin{definition}
    Consider the quaternions $\mathbb{H}$ presented as a complex vector space $\mathbb{H} = \mathbb{C}\oplus \mathbb{C}j.$ The group $\Pin(2)$ is the multiplicative group $U(1)\oplus U(1)j\subset \mathbb{H}.$
\end{definition}

Now, let $X$ be a closed, smooth, oriented 4-manifold with a spin structure. Following Furuta, we will define a map $\mathcal{F}$ similar to $\mathcal{S}$, such that the domain and codomain of $\mathcal{F}$ admit $\Pin(2)$ actions for which $\mathcal{F}$ is equivariant.

Denote by $\mathcal{H}^1$ the space of real-valued 1-forms on $X$ which are harmonic with respect to the metric $g.$ Denote by $\Pi:\Omega^1(\mathbb{R})\to \mathcal{H}^1$ the orthogonal projection coming from the Hodge decomposition.

Suppose the spin$^c$ structure $\mathfrak{s}$ comes from the spin structure on $X$. Then, the determinant bundle $S$ is trivial; choosing a trivialization, we may identify the space of unitary connections $\mathcal{A}(S)$ with the space of purely imaginary 1-forms $\Omega^1(i\mathbb{R})$ on $X.$ We define the map $\mathcal{F}:\Omega^0(W^+\otimes E)\times \Omega^1(\mathbb{R})\to \Omega^0(W^-\otimes E)\times \Omega^0(\mathbb{R})\times \Omega^+(\mathbb{R})\times \mathcal{H}^1$ as follows: $$\mathcal{F}(\Phi, a) = \begin{pmatrix}
    D_{B}^+\Phi + \rho(ia)\Phi \\
    d^*a \\
    (da)^+ +i\rho^{-1}\mu(\Phi, \Phi) \\
    \Pi a
\end{pmatrix}.$$ We may write $$\mathcal{F}(\Phi, a) = D(\Phi, a)+Q(\Phi, a),$$ where: $$D(\Phi, a) = \begin{pmatrix}
    D_{B}^+\Phi \\
    d^*a \\
    (da)^+ \\
    \Pi a \\
\end{pmatrix}, ~  Q(\Phi, a) = \begin{pmatrix}
    \rho(ia)\Phi \\
    0 \\
    i\rho^{-1}\mu(\Phi, \Phi) \\
    0 \\
\end{pmatrix}.$$

Recall that a \emph{real structure} $s$ on a complex vector bundle $V\to X$ is a complex anti-linear endomorphism of $V$ such that $s^2 = 1,$ and a \emph{quaternionic structure} $j$ on a complex vector bundle $V\to X$ is a complex anti-linear endomorphism of $V$ such that $j^2 = -1.$

\begin{proposition}
    Suppose $E$ has a real structure $s$ which is preserved by the connection $B.$ Then, the domain and codomain of $\mathcal{F}$ each admit $\Pin(2)$ actions that are linear over $\mathbb{R}$, and for which $\mathcal{F}$ is equivariant.
\end{proposition}

\begin{remark}
    We wish to explain the intuition for the bizarre but necessary requirement that $E$ have a real structure which is preserved by $B.$ In dimension 4, the spinor bundles $W^\pm$ are in fact \emph{quaternionic}, i.e. they have quaternionic structures $j^\pm$ for which Clifford multiplication is equivariant \cite[Part 2]{Furuta}. We would like to use the quaternionic structure to define a $\Pin(2)$-module structure. Unfortunately, the tensor product of a complex vector bundle with a quaternionic vector bundle does not have a naturally defined quaternionic structure (indeed, the map $j^\pm\otimes \text{id}$ is not well-defined). Hence, we need to assume that our auxiliary bundle $E$ has a real structure $s$; the tensor product $j^\pm\otimes s$ of a quaternionic structure with a real structure is well-defined and quaternionic. The real structure must be parallel with respect to our connection so that the Dirac operator is $\Pin(2)$-equivariant.
\end{remark}

\begin{proof}
    The group $U(1)$ acts by scalar multiplication on sections of $W^\pm\otimes E$, and trivially on the spaces $\Omega^0, \Omega^1,$ and $\Omega^+.$ The map $\mathcal{F}$ is $U(1)$-equivariant with respect to these actions. Therefore, it suffices to construct compatible actions of $j \in \Pin(2)$ on these spaces.

    We first construct a $\Pin(2)$ action on the vector bundles $W^+\otimes E$ and $W^-\otimes E.$ The vector bundles $W^\pm$, being the spinor bundles of a spin structure in dimension 4, have quaternionic structures $j^\pm$ \cite[Part 2]{Furuta}. (Thinking of $W^+$ and $W^-$ as quaternionic vector spaces, these quaternionic structures correspond to right multiplication by $j$.) The tensor product $j^\pm\otimes s$ is therefore a quaternionic structure on $W^\pm\otimes E$ \cite[p. 5]{Crowley-Goette}. We define the action of $j$ on $W^\pm\otimes E$ by the map $j^\pm\otimes s.$ For a section $\sigma\in\Omega^0(E),$ we have $\nabla_B (s\circ\sigma) = s\circ\left(\nabla_B \sigma\right),$ and the untwisted Dirac operator $D^+:\Omega^0(W^+)\to\Omega^0(W^-)$ has $D^+\circ j^+ = j^-\circ D^-,$ hence the twisted Dirac operator $D_B^+$ must be $j$-equivariant.

    On each of the spaces $\Omega^0, \Omega^1, \mathcal{H}^1$ and $\Omega^+,$ we define the action of $j$ to be as $-1$. Since $a\mapsto d^*a+(da)^++\Pi a$ is linear, and $D_B^+\Phi$ is $\Pin(2)$-equivariant, we must have that $D(\Phi, a)$ is $\Pin(2)$-equivariant.

    It remains to show that $Q(\Phi, a)$ is $\Pin(2)$-equivariant. To this end, in each small open set in $X$ we will choose a $\mathbb{C}$-basis of local sections $\sigma_1,\dots, \sigma_r$ of $E$ satisfying $s(\sigma_k) = \sigma_k$ for all $k.$ In these coordinates on $E,$ we may write $\Phi = \sum_{k=1}^r \varphi_k\otimes \sigma_k$ for some sections $\varphi_k\in\Omega^0(W^+).$ Then, $j\Phi = \sum_{k=1}^r (j\varphi_k) \otimes \sigma_k,$ and $$\rho(ia)\Phi = \sum_{k=1}^r (\rho(ia)\varphi_k) \otimes \sigma_k, ~ \mu(\Phi, \Phi) = \sum_{k=1}^r \left(\varphi_k\otimes\varphi_k^* - \frac12|\varphi_k|^2\text{id}\right).$$ 
    
    The $\Pin(2)$-equivariance of $Q$ now follows from the equivariance of the corresponding map $$(\varphi_, a)\mapsto \left(\rho(ia)\varphi, \rho^{-1}\left(\varphi_k\otimes\varphi_k^* - \frac12|\varphi_k|^2\text{id}\right)\right) $$ in ordinary Seiberg-Witten theory established in \cite[Part 2]{Furuta}. 
\end{proof}

Finally, we show the relationship between $\mathcal{F}$ and the multi-monopole moduli space:

\begin{proposition}\label{sw=f}
    Suppose $(X, \mathfrak{s})$ is spin. Then, if the moduli space $\mathcal{M}(\mathfrak{s}, g, B)$ is compact, so is the zero set $\mathcal{F}^{-1}(0)$.
\end{proposition}

\begin{proof}
When $\mathfrak{s}$ comes from a spin structure, the determinant bundle $S$ has a canonical trivialization whose product connection is induced from the metric connection on the spinor bundle.

Let $(\Phi, a)$ be a zero of $\mathcal{F}.$ By our above trivialization of $S,$ the purely imaginary 1-form $ia$ defines a unitary connection $A$ on $S$, defining a map $\mathcal{F}^{-1}(0)\to \mathcal{A}(S)/\mathcal{G}.$ This map is injective, since the Coulomb gauge condition $d^*a = 0$ together with the requirement that $\Pi a = 0$ selects a unique $U(1)$-orbit from the $\mathcal{G}$-orbit of $A$. It is also clear that $(A, \Phi)$ satisfies the multi-spinor Seiberg-Witten equations. Hence, the map $(\Phi, a)\mapsto (A, \Phi)$ is a continuous injection from $\mathcal{F}^{-1}(0)$ to $\mathcal{M}(\mathfrak{s}, g, B)$; since the domain of this injection is closed and the image is compact, we conclude the domain is also compact.
\end{proof}

\section{Finite-Dimensional Approximation of the Multi-Spinor Equations}

Throughout this section, we make the following assumptions: $(X, \mathfrak{s})$ is spin and $E\to X$ has a real structure which is preserved by $B\in\mathscr{A}(E)$.

\subsection{Input from analysis}

Let $V$ be the completion of the domain $\Omega^0(W^+\otimes E)\times \Omega^1$ of $\mathscr{F}$ in the $L^2_4$ topology, and $W$ the completion of the codomain $\Omega^0(W^-\otimes E)\times\Omega^0\times\Omega^+\times \mathcal{H}^1$ of $\mathscr{F}$ in the $L^2_3$ topology. Clearly $\mathscr{F}$ extends to a well-defined map in these completions; moreover $D:V\to W$ is elliptic, so it is Fredholm. By a standard elliptic bootstrapping argument, the zero set $\mathcal{F}^{-1}(0)$ in $L^2_3$ consists solely of smooth zeroes of $\mathcal{F}.$

Let $W'\subset W$ be the $L^2_3$ completion of $\Omega^0(W^-\otimes E)\times\Omega^0\times\Omega^+,$ so that $W = W'\oplus \mathcal{H}^1.$ Following Furuta \cite[p. 282]{Furuta}, we can construct $\Pin(2)$-equivariant norms on $V$ and $W'$, respectively, by: $$||u||_V^2 = \int_X \left(|(D^*D)^2u|^2+|u|^2\right), ~ ||u||_{W'}^2 = \int_X \left(|(DD^*)^{3/2}u|^2+|u|^2\right),$$ and using the $L^2$ norm on $\mathcal{H}^1$ (which induces the $L^2_3$ topology since the elements of $\mathcal{H}^1$ are harmonic) we get a $\Pin(2)$-equivariant norm on all of $W.$

For each $\lambda>0,$ we define $V_\lambda$ to be the span of the eigenvectors of $D^*D$ with eigenvalue $\leq \lambda,$ and likewise define $W_\lambda$ to be the direct sum of $\mathcal{H}^1$ with the span of the eigenvectors of $DD^*$ with eigenvalue $\leq \lambda.$ We let $p_\lambda:W\to W_\lambda$ denote the orthogonal projection. 

We will use the following notation: If $U$ is a normed vector space, for each $R>0,$ we let $S_R(U)$ denote the sphere of radius $R$ in $U$ centered at the origin.

The analytical content of the finite-dimensional approximation method is contained in the following proposition, the main goal of this section:

\begin{proposition}\label{analysis}
    If $\mathcal{M}(\mathfrak{s}, g, B)$ is compact (in particular, if it is empty), then there exists some $R>0$ such that the restriction of $p_\lambda\circ \mathcal{F}$ to $S_R(V_\lambda)$ has no zeroes whenever $\lambda$ is sufficiently large.
\end{proposition}

The proof follows \cite[Lemmas 3.2-3.4]{Furuta} closely, but we write out the details in full below.

\begin{proof}
    By Proposition \ref{sw=f} and our assumption that $\mathcal{M}(\mathfrak{s}, g, B)$ is compact, we have that $\mathcal{F}^{-1}(0)$ is compact. Hence, there exists some $R>0$ such that $\mathcal{F}|S_R(V)$ has no zeroes; we fix a choice of such $R$ for the remainder of this proof.

    \emph{Step 0}: We show that $Q:V\to W$ is a compact operator. The key idea is that $L^2_4(X)$ is a Banach algebra. Hence, if $(\Phi, a)\in V,$ then $Q(\Phi, a)$ is in $L^2_4.$ Since the embedding $L^2_4\hookrightarrow L^2_3$ is compact, we have that $Q:V\to W$ is compact.

    \emph{Step 1}: Suppose $v_k\in V$ is a subsequence converging weakly to $v,$ such that $Dv_k$ converges strongly to $w.$ We show that in fact a subsequence of $v_k$ converges strongly to $v$ and $Dv = w.$ First, since the $v_k$ converge weakly to $v,$ that means $\Pi v_k$ converge strongly to $\Pi v.$ So, by adding the harmonic 1-form $-\Pi v$ to the $v_k$ and to $w$ if necessary, we may suppose that $\Pi v_k = 0$ for all $k.$ Since $L^2_4\hookrightarrow L^2_3$ is compact, up to taking a subsequence if necessary, we conclude $v_k$ converges strongly to $v$ in $L^2_3.$ The elliptic estimate: $$||u||_{L^2_4}\leq C\left(||Du||_{L^2_3}+||u||_{L^2_3}\right)$$ applied to $u = v_k - v_l$ therefore tells us that the sequence $v_k$ is Cauchy in $V,$ so by the uniqueness of weak limits it must converge strongly to $v.$ Since $D:L^2_4\to L^2_3$ is a bounded operator, we must also have that $Dv_k\to Dv$ as $k\to \infty,$ hence $Dv = w.$ 

    \emph{Step 2}: We now show that there exists $\epsilon>0$ so that $||\mathcal{F}v||_W\geq \epsilon$ whenever $v\in S_R(V).$ For, suppose not. Then, there exists a sequence $v_k\in S_R(V)$ such that $||\mathcal{F}v_k||_W\leq \frac1k.$ By the Banach-Alaoglu Theorem, we may pass to a subsequence such that the $v_k$ converge weakly to some $v\in S_R(V)$. Since $Q:V\to W$ is a compact operator, upon taking another subsequence if necessary, we may assume that $Qv_k$ converges strongly to $Qv.$ Since $$||Dv_k+Qv_k||_W\leq \frac1k$$ for each $k,$ we conclude that $Dv_k$ converges strongly to $-Qv.$ By Step 1, we conclude that $Dv = - Qv,$ which contradicts that there are no zeroes of $D+Q$ with norm $R.$

    \emph{Step 3}: We show that for any fixed $w\in W,$ we have $$\lim_{\lambda\to\infty} (1-p_\lambda)w = 0.$$ Let $\lambda_k,$ for $k\in\{1,2,\dots\}$ denote the spectrum of the self-adjoint elliptic operator $DD^*$, as a non-decreasing sequence where each eigenvalue is repeated according to its multiplicity. By the Spectral Theorem, we may write $w = \sum_{k} a_k u_k,$ where $u_k$ are eigenfunctions of the operator $DD^*$ with eigenvalue $\lambda_k,$ scaled to have $||u_k||_{L^2} = 1,$ and the series converges in $L^2.$ We may further assume that if $k\neq j$ then $\langle u_k, u_j\rangle_{L^2} = 0$. Since $w \in W,$ we have that $(DD^*)^{3/2}w\in L^2.$ Since $$\langle (DD^*)^{3/2}w, u_k\rangle = \langle w, (DD^*)^{3/2}u_k\rangle = \lambda_k^{3/2}a_k,$$ we must have that in $L^2$, we have $$(DD^*)^{3/2}w = \sum_{k} \lambda_k^{3/2}a_k u_k.$$ Hence, $$||w||_W^2 = ||w||_{L^2}^2 + ||(DD^*)^{3/2}w||_{L^2}^2 = \sum_k\left(1+\lambda_k^3\right)|a_k|^2.$$

    Each $u_k$ is an eigenfunction of an elliptic operator, so it is smooth, so a fortiori it is in the space $L^2_3.$ We may compute the inner product of these eigenfunctions in $W$ as follows: $$\langle u_k, u_j\rangle_{W} = \langle u_k, u_j\rangle_{L^2} + \langle (DD^*)^{3/2}u_k, (DD^*)^{3/2}u_j\rangle_{L^2} = \left(1+ \lambda_k^{3/2}\lambda_j^{3/2}\right)\langle u_k, u_j\rangle_{L^2},$$ so we find if $k\neq j$ then $u_k$ and $u_j$ are also orthogonal in $W$. Hence, $$||p_\lambda w||_W^2 = \sum_{\lambda_k\leq\lambda} \left(1+\lambda_k^3\right)|a_k|^2$$ and $p_\lambda w$ is orthogonal to $(1-p_\lambda)w.$ So, the Pythagorean Theorem gives us: $$||(1-p_\lambda)w||_W^2 = ||w||_W^2  - ||p_\lambda w||_W^2 = \sum_{\lambda_k>\lambda} \left(1+\lambda_k^3\right)|a_k|^2,$$ which converges to 0 as $\lambda\to \infty$ since it is the tail of a convergent series.

    \emph{Step 4}: We show that for all $\lambda$ sufficiently large, the following is true for all $v\in S_R(V)$: $$||(1-p_\lambda)Qv||<\epsilon. $$ For, if not, there exists a sequence $v_k\in S_R(V)$ with: $$||(1-p_k)Qv_k||\geq\epsilon$$ for each positive integer $k.$ Upon choosing a subsequence if necessary, we may assume $v_k$ converges weakly to some $v\in S_R(V)$ and $Qv_k$ converges strongly to $Qv.$ Therefore, \begin{equation*}
        \begin{split}
            \epsilon &< ||(1-p_k)Qv_k|| \\
            &\leq ||(1-p_k)Qv||+||(1-p_k)\left(Qv_k - Qv\right)|| \\
            &\leq  ||(1-p_k)Qv||+||Qv_k - Qv|| \\
        \end{split}
    \end{equation*} As $k\to \infty,$ the quantity $||(1-p_k)w||$ approaches zero for any fixed $w\in W$ by the result of Step 3. Since $||Qv_k - Qv||$ also approaches zero as $k\to \infty,$ we arrive at a contradiction.

    \emph{Step 5}: Choose $\lambda$ sufficiently large in the sense of Step 4. Then, if $v\in S_R(V_\lambda),$ clearly $Dv \in W_\lambda$, so we combines Steps 2 and 4 to conclude: $$||(p_\lambda\circ\mathcal{F})(v)|| =  ||Dv+p_\lambda Qv|| \geq ||Dv+ Qv|| - ||(1-p_\lambda)Qv||>\epsilon - \epsilon = 0.$$ This completes the proof.
\end{proof}

\subsection{$\Pin(2)$-Equivariant Index Theory}

We fix $R$ and $\lambda$ from Proposition \ref{analysis}. Our next task is to compute the decomposition of the finite-dimensional $\Pin(2)$-modules $V_\lambda$ and $W_\lambda$ into irreducible $\Pin(2)$-modules. We will see that they are direct sums of the following two Lie group representations: \begin{enumerate}
    \item $\widetilde{\mathbb{R}}$ is the 1-dimensional real representation of $\Pin(2)$ uniquely specified by the fact that $j$ acts by $-1.$
    \item $\mathbb{H}$ is the 4-dimensional real representation of $\Pin(2)$ given by left-multiplication on the quaternions, considered as a real vector space.
\end{enumerate}

We also use the following standard notations for the well-known topological invariants of the $4$-manifold $X$:
\begin{itemize}
    \item $b_1(X)$ denotes the rank of $H^1(X;\mathbb{R})$,
    \item $b_+(X)$ is the rank of a maximal positive definition subspace of $H^2(X;\mathbb R)$,
    \item $\sigma(X)$ refers to the signature of $X$.
\end{itemize}
When context is clear, we sometimes drop the dependency on the $4$-manifold and simply write these invariants as $b_1, b_+$, and $\sigma$.

\begin{proposition}
    Let $\kappa = -\frac{r\sigma(X)}{16}-\frac12c_2(E).$ There exist integers $m, n$ such that: $$V_\lambda = \mathbb{H}^{\kappa+m}\oplus \widetilde{\mathbb{R}}^n, ~W_\lambda = \mathbb{H}^{m}\oplus \widetilde{\mathbb{R}}^{b_+(X)+b_1(X)+1+n}$$
\end{proposition}

\begin{proof}
    Decompose the operator $D$ as $D_1+D_2,$ where $D_1 = D_B^+:\Omega^0(W^+\otimes E)\to\Omega^0(W^-\otimes E)$ and $D_2 = d^+\oplus d^*\oplus\Pi:\Omega^1\to \Omega^+\oplus\Omega^0\oplus \mathcal{H}^1.$

    We first decompose the spaces $V_\lambda, W_\lambda$ into orthogonal direct sums $V_\lambda = V_\lambda^1\oplus V_\lambda^2$ and $W_\lambda = W_\lambda^1\oplus W_\lambda^2$ as follows: for $i= 1,2,$ we let $V^i, W^i$ be the domain and codomain, respectively, of the map $D_i,$ completed in the norm on $V, W$ respectively. We then let the spaces $V_\lambda^i$ be the intersections $V_\lambda^i := V_\lambda\cap V^i,$ and likewise for the spaces $W_\lambda^i.$ 

    Following \cite[Lemma 4.1]{Furuta}, we note that there are finitely many points $p_1,\dots, p_\ell$ such that the evaluation map $$\bigoplus_{i=1}^\ell\text{ev}_{p_i}:V_\lambda^1\to \bigoplus_{i=1}^\ell (W^+\otimes E)_{p_i}$$ is an injection of $\Pin(2)$-modules. Since each fiber $(W^+\otimes E)_{p_i}$ is isomorphic to $\mathbb{H}^r$ as a $\Pin(2)$-module, the right-hand side is isomorphic to $\mathbb{H}^{\ell r}$ as a $\Pin(2)$-module. Since $\mathbb{H}$ is an irreducible representation, the left-hand side $V_\lambda^1$ must also be isomorphic to a direct sum of copies of $\mathbb{H}.$ Likewise, $W_\lambda^1$ is also isomorphic to a direct sum of copies of $\mathbb{H}.$ The difference between the multiplicity of $\mathbb{H}$ in $V_\lambda^1$ and $W_\lambda^1$ is therefore given by half the difference in their complex dimensions. But this is precisely $\frac12 \ind_\mathbb{C} D_B^+.$

    An identical argument with $V_\lambda^2$ and $W_\lambda^2$ shows that these $\Pin(2)$-representations are each isomorphic to a direct sum of  $\widetilde{\mathbb{R}}$, and that the difference between the multiplicity of $\widetilde{\mathbb{R}}$ in $V_\lambda^2$ and that in $W_\lambda^2$ is given by $\ind_\mathbb{R} (d^*+d^+:\Omega^1\to\Omega^0\oplus\Omega^+)-\dim \mathcal{H}^1.$ By the Atiyah-Singer Index Theorem \cite{AS3}, the former quantity is $$\frac12\ind_\mathbb{C} D_B^+ =-\frac{r\sigma(X)}{16}-\frac12c_2(E),$$ while the latter quantity is $-b^++b_1-1-b_1 = -b^+-1.$
\end{proof}

In fact, the image of $p_\lambda\circ\mathcal{F}$ is in the orthogonal complement of the constant function $1\in\Omega^0$, in $W.$ Since the restriction of the $\Pin(2)$ action to $\mathbb{R}\cdot 1$ is a copy of $\widetilde{\mathbb{R}}$, we may replace $W_\lambda$ with $$\overline{W}_\lambda:= (\mathbb{R}\cdot 1)^\perp = \mathbb{H}^{\kappa}\oplus \widetilde{\mathbb{R}}^{b_+(X)+n}.$$ So, putting the results of this section together we arrive at the following: 

\begin{proposition}\label{mswsetup}
    Let $\kappa = -\frac{r\sigma(X)}{16}-\frac12c_2(E).$ Suppose there exists a metric $g$ on $X$ and a $SU(r)$-compatible connection $B$ on $E$ such that $\mathcal{M}(g, B)$ is compact. Then, there exist integers $m, n,$ a real number $R>0,$ and a $\Pin(2)$-equivariant map $$p_\lambda\circ\mathcal{F}:\mathbb{H}^{\kappa+m}\oplus \widetilde{\mathbb{R}}^n\to \mathbb{H}^{m}\oplus \widetilde{\mathbb{R}}^{b_+(X)+n}$$ that has no zeroes on $S_R(\mathbb{H}^{\kappa+m}\oplus \widetilde{\mathbb{R}}^n).$
\end{proposition}

\subsection{Generalized Furuta Theorem}

Furuta uses the Adams operations in equivariant $K$-theory to prove the following proposition (see also Bryan's beautiful reproof \cite[Section 3.4]{Bryan}, which we include in Appendix A):

\begin{proposition}[{\cite[Proposition 5.1]{Furuta}}]\label{ktheory}
    Suppose there exist nonnegative $m, n,$ positive integers $\beta, \kappa,$ a real number $R>0,$ and a $\Pin(2)$-equivariant map $$f:\mathbb{H}^{\kappa+m}\oplus \widetilde{\mathbb{R}}^n\to \mathbb{H}^{m}\oplus \widetilde{\mathbb{R}}^{\beta+n}$$ that has no zeroes on $S_R(\mathbb{H}^{\kappa+m}\oplus \widetilde{\mathbb{R}}^n).$ Then, $$2\kappa+1\leq \beta.$$
\end{proposition}

Combining the above proposition with Proposition \ref{mswsetup} (and introducing an extraneous factor of $\widetilde{\mathbb{R}}$ on the right-hand side if $b_+(X) = 0$), we arrive at:

\begin{theorem}\label{Furuta}
    If $\mathcal{M}(g, B)$ is compact (in particular, if it is empty), then: $$-\frac{r\sigma(X)}{16}-\frac12 c_2(E) \leq \max\left(0, \frac{b_+(X)-1}{2}\right).$$ So, if $b_+(X)+b_1(X)\geq 1$, $$-\frac{r+4}4\sigma(X)-2c_2(E) \leq b_2(X) - 2.$$
\end{theorem}

\section{Applications to Algebraic Geometry}

In this section, we will show how to combine Theorem \ref{Furuta} with algebro-geometric interpretations of the multi-spinor Seiberg-Witten equations to produce the main theorems of this paper.

\begin{remark}
    Our lower bounds on the second Chern classes of good vector bundles do not require the bundles to be stable. However, the geometric interpretation of the ``good'' condition in terms of the smoothness of moduli spaces, both in gauge theory and in algebraic geometry, becomes much clearer when we assume stability. Hence, we will frequently restrict to the stable case and compare our bound to other bounds that do require stability.
\end{remark}

\subsection{Holomorphic structures and connections}

Fix a K\"ahler surface $(X, \omega)$ and an $SU(r)$-bundle $E\to X.$ A connection $B\in\mathscr{A}(E)$ is \emph{(1,1)} if $F_B^{2,0} = 0$. In this case it defines a holomorphic structure on $E$, and we denote the corresponding holomorphic bundle by $E_B.$ A (1,1)-connection $B$, or its corresponding holomorphic vector bundle $E_B$, is: \begin{enumerate}
    \item \emph{good} if $H^0(\End_0 E_B\otimes K_X)=0$.
    \item \emph{semi-weakly good} if every non-trivial section $s\in H^0(\End E_B\otimes K_X)$ has rank $\geq 2$ at some point.
    \item \emph{weakly good} if every non-trivial section $s\in H^0(\End_0 E_B\otimes K_X)$ has rank $\geq 2$ at some point.
\end{enumerate} 

Note that these definitions make sense for a holomorphic bundle without a prescribed connection as well, and that a holomorphic bundle $\mathscr{E}\to X$ is good (resp. semi-weakly good, weakly good) if and only if it admits a hermitian metric whose Chern connection is good (resp. semi-weakly good, weakly good.)

\begin{remark}
    The definition of \emph{good} here is that $H^0(\End_0 E_B\otimes K_X)=0$, while the definition of \emph{good} we gave in Definition \ref{good} is that $H^2(\End_0 E_B)=0$. Since $\End_0 E_B$ is holomorphically isomorphic to its dual, Serre duality identifies $H^2(\End_0 E_B)$ with the dual space of $H^0(\End_0 E_B\otimes K_X)$ so the definitions agree.
\end{remark}

\begin{remark}
    If the rank of $E$ is $\geq 2,$ then (1) implies (2) implies (3). That (2)$\implies$(3) follows immediately from the direct sum decomposition $\End E_B = \End_0 E_B \oplus \mathbb{C}\cdot\text{id}$. To see (1)$\implies$(2), note that if $H^0(\End_0 E_B\otimes K_X)=0,$ the only sections of $H^0(\End E_B\otimes K_X)=0$ are products of the identity with a section of $K_X.$ Hence, if the rank of $E$ is at least 2, then any nonzero section of $\End E_B\otimes K_X$ must have rank at least 2 at some point.
\end{remark}

\begin{remark}
    We give a gauge-theoretic interpretation of a ``good'' vector bundle. A stable holomorphic vector bundle $\mathscr{E}$ admits an anti-self-dual Yang-Mills connection $B$ (see for instance \cite[Theorem 6.1.5]{DK}.) This is a (1,1)-connection and and if it is also irreducible the second cohomology of the deformation complex of the anti-self-duality equations at $B$ may be identified with $H^2(\End_0 E_B)$ \cite[p. 240]{DK}. The second cohomology is the obstruction space in the Kuranishi model of instanton moduli space, so we may think of good connections $B$ as ``unobstructed.'' In particular, instanton moduli space is smooth in a neighborhood of $[B]$ whenever $B$ is irreducible and $H^2(\End_0 E_B) = 0$.
\end{remark}

\begin{theorem}\label{compactness}
    Let $(X, \omega)$ be a K\"ahler surface with metric $g$ and spin structure determined by a choice of square root of the canonical bundle $K^{1/2}$. If $B$ is a weakly good (1,1)-connection on $E$, then $\mathcal{M}(g, B)$ is compact.
\end{theorem}

\begin{proof}

    By \cite[Equation 2, p. 7]{Thakar26}, the multi-spinor Seiberg-Witten equations reduce to the following equations, for a unitary connection $A$ on the bundle $K^{1/2}$, a section $\alpha\in\Omega^0(E\otimes L),$ and a section $\beta\in\Omega^{0,2}(E\otimes L)$:

    \begin{equation}\label{holomsw}\begin{cases}
    \delbar_{AB}\alpha &= -\delbar_{AB}^*\beta \\
    F_A^{2,0} &= \tr(\alpha\otimes\beta^*) \\
    -\sqrt{-1}\Lambda F_A &= \frac14\left(|\alpha|^2-|\beta|^2\right) - \frac12\sqrt{-1}\Lambda F_{\nabla_0}.\\
\end{cases}\end{equation}

    By our assumption that $F_B^{2,0} = 0,$ \cite[Lemma 3.1]{Thakar26} shows any solution has $$\delbar_{AB}\alpha = \delbar_{AB}^*\beta = F_A^{2,0} = \tr(\alpha\otimes\beta^*) = 0.$$ Denote by $\mathscr{E}$ the holomorphic structure $E_B.$ Then, $\alpha\otimes\beta^*\in H^0_X(\End \mathscr{E}\otimes K_X)$ is traceless, so it defines a section of $H^0_X(\End_0 \mathscr{E}\otimes K_X),$ Because $\mathscr{E}$ is weakly good, all elements of this space with rank $\leq 1$ pointwise are identically zero. Hence, $\alpha\otimes\beta^* = 0,$ so at each point $x\in X,$ one of $\alpha(x)$ or $\beta(x)$ must equal 0. Since they are both holomorphic, one must vanish identically by unique continuation. 
    
    By Chern-Weil Theory, $\left[\frac{i}{2\pi}F_A\right] = c_1(K^{1/2})\in H^2(X; \mathbb{R}).$ Hence, the integral of the third holomorphic Seiberg-Witten equation then gives us: $$0 =\pi\left(2c_1(K_X^{1/2})-c_1(K_X)\right)\cdot[\omega] = ||\alpha||_{L^2}^2 - ||\beta||_{L^2}^2.$$ Since at least one of $\alpha$ or $\beta$ must vanish, we in fact have that $\alpha = \beta = 0.$ This implies we may identify the space of solutions with the space of connections $A$ satisfying the equations $$ F_A^{2,0} = 0, ~ \Lambda F_A = \frac12\Lambda F_{\nabla_0}$$ modulo unitary gauge, which, as is well-known, may be identified with the torus $H^1(X; \mathbb{R}/\mathbb{Z})$. This latter space is evidently compact.
\end{proof}

\begin{remark}
    This theorem is a variant of \cite[Theorem 4.3]{Thakar26}. The theorem \cite[Theorem 4.3]{Thakar26} tells us that the moduli space of \emph{perturbed} solutions converges when the perturbation is sufficiently large. Here we only consider unperturbed solutions.
\end{remark}

\subsection{An Example of a Weakly Good Bundle}

The compactness theorem holds for weakly good (1,1)-connections, which are a generalization of good (1,1)-connections. We see in the following lemma an example of weakly good bundles which are not good:

\begin{lemma}\label{directsum}
    Let $X$ be a K\"ahler surface and $\mathscr{E}\to X$ a semi-weakly good holomorphic bundle of rank $r\geq 2$. Then, for each positive integer $\ell,$ the bundle $\zeta_\ell\to X$ defined as the connected sum of $\ell$ copies of $\mathscr{E}$ is semi-weakly good.
\end{lemma}

\begin{proof}
    Note that $\End \zeta_\ell$ is isomorphic to the direct sum of $\ell^2$ copies of $\End \mathscr{E},$ corresponding to the decomposition of an endomorphism of $\zeta_\ell$ into block matrices corresponding to each $\mathscr{E}$ factor. Since $\mathscr{E}$ is semi-weakly good, every nonzero section of $\End \mathscr{E}\otimes K$ has rank at least $r\geq 2$ at some point. Therefore, every nonzero section of $\End \zeta_\ell$ also has rank at least 2 at some point.
\end{proof}

\begin{remark}
    Note that if $\ell\geq 1$ then $\zeta_\ell$ is obstructed whenever there exists a nonzero section of the canonical bundle $\sigma\in H^0(K)$. To see this, label the copies of $\mathscr{E}$ in $\zeta_\ell$ as $\mathscr{E}_1, \mathscr{E}_2,\dots, \mathscr{E}_\ell.$ Then, $$\text{id}_{\mathscr{E}_1}\otimes\sigma - \text{id}_{\mathscr{E}_2}\otimes\sigma$$ is a nonzero section of $(\End_0 \zeta_\ell) \otimes K.$
\end{remark}

\subsection{Real bundles on K\"ahler surfaces}
On a K\"ahler surface $X$, $b_+(X) = 1+2h^2(\mathscr{O}_X)\geq 1,$ so we may combine our compactness result, Theorem \ref{compactness}, with the result of our finite-dimensional approximation, Theorem \ref{Furuta}, to get:

\begin{corollary}\label{realkahler}
    Let $X$ be a spin K\"ahler surface and let $\mathscr{E}$ be a weakly good holomorphic bundle on $X$ with trivial determinant, and with a real structure and compatible holomorphic connection. Then, $$\frac{\rank \mathscr{E} +4}{4}(-\sigma(X))-2c_2(\mathscr{E}) \leq b_2(X)- 2.$$
\end{corollary}

With a little more work, we get the following stronger corollary:

\begin{corollary}\label{realstable}
    Let $X$ be a spin K\"ahler surface and let $\mathscr{E}$ be a semi-weakly good holomorphic bundle on $X$ with trivial determinant, and with a real structure and compatible holomorphic connection. Then, $$c_2(\mathscr{E})\geq \rank \mathscr{E}\frac{-\sigma(X)}{8} = \rank \mathscr{E}\left(\chi(\mathscr{O}_X) - \frac18c_1^2(X)\right).$$
\end{corollary}

\begin{proof}
    Define $\zeta_\ell\to X$ to be the direct sum of $\ell$ copies of $\mathscr{E},$ for some positive integer $\ell.$ Since a direct sum of real connections is real, $\zeta_\ell$ has a real structure with compatible connection. By Lemma \ref{directsum}, $\zeta_\ell$ is weakly good. Hence, Corollary \ref{realkahler} tells us: $$\frac{\rank \zeta_\ell +4}{4}(-\sigma(X))-2c_2(\zeta_\ell) \leq b_2(X)+2b_1(X) - 2.$$ Using the relations $\rank \zeta_\ell = \ell\rank \mathscr{E}$ and $c_2(\zeta_\ell) = \ell c_2(\mathscr{E}),$ we conclude: $$\frac{\ell\rank \mathscr{E} +4}{4}(-\sigma(X))-2\ell c_2(\mathscr{E}) \leq b_2(X)+2b_1(X) - 2.$$ Rearranging, we get: $$c_2(\mathscr{E}) \geq \frac{\ell\rank \mathscr{E} +4}{8\ell}(-\sigma(X))- \frac{b_2(X) +2b_1(X)+ 2}{2\ell}.$$ Since this is true for every positive integer $\ell$, we may take the limit as $\ell\to\infty$ to conclude that $$c_2(\mathscr{E})\geq \rank \mathscr{E}\frac{-\sigma(X)}{8},$$ and that suffices for the proof.
\end{proof}

Another source of real vector bundles is from vector bundles that are holomorphically isomorphic to their dual bundle; this results in the following theorem:

\begin{theorem}\label{selfdual}
    If $X$ is a spin K\"ahler surface and $\mathscr{E}\to X$ is a good holomorphic vector bundle of rank $r$ which is holomorphically isomorphic to its dual bundle $\mathscr{E}^*\to X,$ then: $$c_2(\mathscr{E}) \geq -r\frac{\sigma(X)}{8} = r\left(\chi(\mathscr{O}_X) - \frac18 c_1^2(X)\right)$$
\end{theorem}

\begin{proof}
    the bundle $\mathscr{E}\oplus\mathscr{E}^*\to X$ is isomorphic to the bundle $\mathscr{E}\oplus\mathscr{E},$ so it is semi-weakly good by Lemma \ref{directsum}. By Theorem \ref{sp}, it suffices to show that the bundle $\mathscr{E}\oplus\mathscr{E}^*\to X$ has a parallel real structure.

    Since $\mathscr{E}$ is holomorphically isomorphic to its dual, we must have $c_1(\mathscr{E}) = -c_1(\mathscr{E}^*) = 0,$ so upon taking the tensor product of $\mathscr{E}$ with a trivial holomorphic line bundle if necessary, we may assume that $\mathscr{E}$ has trivial determinant. Let $E$ denote the underlying vector bundle; it therefore has $SU(r)$ structure group and a compatible Chern connection $B$. Let $\tau:E\to E^*$ denote the holomorphic isomorphism from $E$ to $E^*.$ The map $s:E\oplus E^*\to E\oplus E^*$ given by $s(\alpha, \beta) = (\beta^*, \alpha^*)$ is complex-antilinear with $s^2 =1,$ and it is clearly parallel, sufficing for our claim.
\end{proof}

Theorem \ref{selfdual} has a curious implication, which we now describe. Let $(X, \omega)$ be a K\"ahler surface and let $M_X^{good}(r, c)$ denote the subspace of moduli space of $\omega$-stable bundles $\mathscr{E}\to X$ with rank $r$, trivial determinant, and $c_2(\mathscr{E}) = c$ consisting of only those bundles which are good. For any good $\omega$-stable rank $r$ bundle $\mathscr{E}\to X$ with trivial determinant, the dual bundle $\mathscr{E}^*$ is $\omega$-stable, it has $c_2(\mathscr{E}^*) = c_2(\mathscr{E})$ since $c_1(\mathscr{E}) = 0,$ and it is good since $\End_0\mathscr{E}\cong \End_0\mathscr{E}^*.$ Hence, there is a well-defined involution $\tau:M_X^{good}(r, c)\to M_X^{good}(r, c)$ given by $[\mathscr{E}]\mapsto [\mathscr{E}^*]$.

\begin{corollary}
    Let $X$ be a spin K\"ahler surface. If $c< r\left(\chi(\mathscr{O}_X) - \frac18 c_1^2(X)\right),$ then the involution $\tau:M_X^{good}(r, c)\to M_X^{good}(r, c)$ is free.
\end{corollary}

\begin{proof}
    If the involution $\tau$ fixes an equivalence class $[\mathscr{E}]\in M_X^{good}(r, c)$, that implies $\mathscr{E}$ is holomorphically isomorphic to its dual bundle. Theorem \ref{selfdual} therefore implies that $$c = c_2(\mathscr{E})\geq r\left(\chi(\mathscr{O}_X) - \frac18 c_1^2(X)\right),$$ contradicting our assumption.
\end{proof}

\subsection{Quaternionic bundles on K\"ahler surfaces}

In this part, we will extract characteristic class inequalities for good holomorphic bundles not necessarily having a real structure. Our strategy is to build real bundles out of the complex bundles, using the following lemma:

\begin{lemma}\label{complex-to-real}
    If $E\to X$ is a quaternionic bundle with compatible connection $B,$ then the bundle $E\oplus E \to X$ has a real structure compatible with the connection $B\oplus B$.
\end{lemma}

\begin{proof}
    Let $j:E\to E$ be the parallel quaternionic structure, so that $j$ is a complex anti-linear map with $j^2 =-1$. The map $(\alpha, \beta)\mapsto (-j\beta, j\alpha)$ on $E\oplus E$ is then a complex anti-linear involution which is parallel with respect to the connection $B\oplus B.$
\end{proof}

\begin{theorem}\label{sp}
    Let $X$ be a spin K\"ahler surface and let $\mathscr{E}$ be a semi-weakly good holomorphic bundle on $X$ with a quaternionic structure and connection compatible with both the quaternionic structure and the holomorphic structure. Then, $$c_2(\mathscr{E})\geq \rank \mathscr{E}\frac{-\sigma(X)}{8} = \rank \mathscr{E}\left(\chi(\mathscr{O}_X) - \frac18 c_1^2(X)\right).$$
\end{theorem}

\begin{proof}
    By Lemma \ref{complex-to-real}, the bundle $\mathscr{E}\oplus \mathscr{E}$ has a real structure with compatible holomorphic connection, and by Lemma \ref{directsum}, it is semi-weakly good. Since $Sp(n) \subseteq SU(2n),$ the bundle $\mathscr{E}$ has $c_1(\mathscr{E}) = 0,$ so $c_2(\mathscr{E}\oplus \mathscr{E}) = 2c_2(\mathscr{E})$. The desired inequality therefore follows immediately from Theorem \ref{realstable} together with the relations $\rank \mathscr{E}\oplus \mathscr{E} = 2 \rank \mathscr{E}$ and $c_2(\mathscr{E}\oplus \mathscr{E}) = 2c_2(\mathscr{E})$.
\end{proof}

As a consequence of the exceptional isomorphism $SU(2)\cong Sp(1)$, any $SU(2)$ bundle with compatible connection automatically has a parallel quaternionic structure. To see this more concretely, if $E$ is a rank 2 bundle with a fixed trivialization of the determinant $\Lambda^2E\cong \underline{\mathbb{C}}$, the map $\tau:E\to E^*$ given by $\tau(v)(w) := v\wedge w$ is a complex-linear isomorphism. Composing the standard anti-holomorphic isomorphism $E\to E^*$ given by the hermitian metric with $\tau^{-1}$ gives us a complex anti-linear map $j:E\to E$ which is evidently parallel. From the local coordinate expression $j(a, b) = (-\overline b, \overline a),$ we can immediately verify that $j^2 = -1$ so it is indeed a quaternionic structure. We get the following corollary:

\begin{corollary}\label{su(2)}
    Let $X$ be a spin K\"ahler surface and let $\mathscr{E}$ be a good bundle on $X$ with rank $2$ and trivial determinant. Then, $$c_2(\mathscr{E})\geq \frac{-\sigma(X)}{4} = 2\chi(\mathscr{O}_X) - \frac14 c_1^2(X).$$ In other words, $$n_2(X) \geq 2\chi(\mathscr{O}_X) - \frac14 c_1^2(X).$$
\end{corollary}

\begin{remark}
    We can also see Corollary \ref{su(2)} as a corollary of Theorem \ref{selfdual}, since the map $\tau:E\to E^*$ described above is in fact a holomorphic isomorphism.
\end{remark}

\subsection{Vector bundles on K3 surfaces}

In this subsection, we consider the case when $X$ is a K3 surface, i.e. it has trivial canonical bundle and $h^1(\mathscr{O}_X) = 0.$ Since $K_X$ is trivial, it is in particular even, so $X$ is spin. The inequalities we will derive with our methods below, while already known, are sharp. 

What makes this case notable is that the trivial canonical bundle in combination with Serre duality can simplify many sheaf cohomology computations. We gather these remarkable conclusions in the following lemma:

\begin{lemma}\label{k3simple}
    Let $(X,\omega)$ be a K\"{a}hler surface with trivial canonical bundle, and let $\mathscr{E}\to X$ and $\mathscr{F}\to X$ be $\omega$-stable vector bundles with rank $r\geq 2$ and with $c_1(\mathscr{E})=c_1(\mathscr{F}) = 0.$ \begin{enumerate}
        \item $\mathscr{E}$ is good.
        \item All nontrivial sections in $H^0(\Hom(\mathscr{E}, \mathscr{F}))$ have pointwise rank $\geq 2.$
        \item The bundle $\mathscr{E}\oplus \mathscr{E}^*$ is semi-weakly good.
    \end{enumerate}
\end{lemma}

\begin{proof}
    A main idea in the proofs of these lemmas is the statement that stable bundles are simple, i.e. that if $\mathscr{E}\to X$ and $\mathscr{F}\to X$ are stable vector bundles with the same rank and degree, then $H^0(\Hom(\mathscr{E}, \mathscr{F}))$ is 0 if they are not isomorphic, and it consists entirely of constant multiples of the identity if they are isomorphic.
    \begin{enumerate}
        \item Since $K_X$ is trivial, $H^0(\End_0\mathscr{E}\otimes K_X)= H^0(\End_0\mathscr{E}),$ and this latter space is zero by the simplicity of $\mathscr{E}.$
        \item This follows immediately from the simplicity of $\mathscr{E}, \mathscr{F}$ together with the observation that $r\geq 2.$
        \item If $\mathscr{E}$ is stable of $c_1(\mathscr{E}) = 0,$ then its dual bundle $\mathscr{E}^*$ is stable and also has $c_1(\mathscr{E}^*) = 0$. Hence, by item (2), $H^0(\Hom(\mathscr{E}, \mathscr{E}^*)) = H^0(\Hom(\mathscr{E}, \mathscr{E}^*)\otimes K_X)$ consists only of the zero section, together with elements with pointwise rank $\geq 2.$
    \end{enumerate}
    This suffices for the proof.
\end{proof}

\begin{theorem}
    Let $X$ be a K3 surface and $\mathscr{E}\to X$ a stable vector bundle with trivial determinant. Then, $$c_2(\mathscr{E})\geq 2\rank \mathscr{E}.$$ In other words, $$n_r(X) \geq 2r.$$
\end{theorem}

\begin{proof}
    Let $E\to X$ be the corresponding $SU(r)$-bundle and $B\in \mathscr{A}(E)$ the corresponding anti-self-dual connection. The bundle $E\oplus E^*$, with connection $B\oplus B^*,$ has a parallel real structure $(\alpha, \beta)\mapsto (\beta^*, \alpha^*)$. It corresponds to the holomorphic bundle $\mathscr{E}\oplus \mathscr{E}^*$, which is semi-weakly good by Lemma \ref{k3simple}. Hence, Theorem \ref{realstable} tells us $c_2(\mathscr{E}\oplus \mathscr{E}^*)\geq \rank(\mathscr{E}\oplus \mathscr{E}^*)\frac{-\sigma(X)}{8}.$ The result follows from the additivity of rank and of second Chern class when $c_1 = 0,$ and the well-known fact that the signature of a K3 surface is $-16.$
\end{proof}

This inequality has in fact been proven by much simpler means, by T. Nakashima, who in addition proved it was sharp for algebraic K3 surfaces.

\begin{theorem}[{\cite[Theorem 2.3]{Nakashima}}]
    Let $(X, H)$ be a polarized K3 surface. Then, $n_r(X, H) = 2r.$
\end{theorem}

\begin{proof}[Proof of lower bound on $n_r$]
    Suppose $\mathscr{E}\to X$ is a good (equivalently, stable) bundle with rank $r$, $c_1(\mathscr{E}) = 0,$ and $c_2(\mathscr{E}) = c$. Then, $\mathscr{E}^*$ is also stable with $c_1(\mathscr{E}^*) = 0.$ So, by Serre duality and the fact that a degree 0 stable vector bundle of rank $\geq 2$ has no nonzero sections: $$H^0(\mathscr{E}) = 0, ~H^2(\mathscr{E}) = H^0(\mathscr{E^*}\otimes K) = H^0(\mathscr{E}^*) = 0.$$ By the Hirzebruch-Riemann-Roch Theorem: $$-c + 2r = \chi(\mathscr{E}) = h^2(\mathscr{E}) - h^1(\mathscr{E}) + h^0(\mathscr{E}) = -h^1(\mathscr{E})\leq 0, $$ as desired.
\end{proof}

\subsection{The Gieseker-Friedman Bounds and Elliptic Surfaces}

In the remainder of this section, we will take stock of Corollary \ref{su(2)} in the context of existing bounds on the invariant $n_2$, and in more generality $n_r$, in the algebraic geometry literature.

Numerous authors have found upper bounds on $n_2$; this amounts to constructing good vector bundles on $X$ and computing their second Chern classes. A very strong such bound was found by Gieseker for algebraic surfaces (where we use the notation $(X, H)$ for a polarized algebraic surface synonymously with $(X, \omega),$ where $\omega$ is the associated K\"ahler form):

\begin{theorem}[{\cite[Theorems 1.1, 1.2]{Gieseker}}]
    Let $(X, H)$ be a polarized smooth surface. \begin{enumerate}
        \item $n_2(X, H)\leq \begin{cases} 2p_g(X)+4 & \text{if}~p_g(X)~\text{is even} \\ 2p_g(X)+2 & \text{if}~p_g(X)~\text{is odd.}\end{cases} $
        \item If $p_g(X)>1000,$ then $n_2(X, H) \leq \frac32p_g(X) + 6.$
    \end{enumerate}
\end{theorem}

For elliptic surfaces this was strengthened by Friedman: 

\begin{theorem}[{\cite[Theorem 0.4]{Friedman-elliptic}}]
    Let $(X, H)$ be a polarized elliptic surface with $h^1(\mathscr{O}_X) = 0.$ We have: $$n_2(X, H)\leq 2p_g(X) + 2.$$
\end{theorem}

For a surface $X$, our lower bound on $n_2$ in Corollary \ref{su(2)} may be re-written as: $$n_2(X)\geq 2p_g+2-2h^1(\mathscr{O}_X) - \frac14c_1^2(X).$$ By comparing Gieseker's bound to ours, we find that it is sharp only when $2h^1(\mathscr{O}_X)+\frac14 c_1^2(X) = 0.$ 

\begin{lemma}
    If $X$ is a spin K\"ahler surface that satisfies $2h^1(\mathscr{O}_X)+\frac14 c_1^2(X) = 0,$ then it is either a K3 surface or it is minimal and elliptic.
\end{lemma}

\begin{proof}
    Note that $$-\frac{\sigma(X)}4 = 2p_g+2-2h^1(\mathscr{O}_X) - \frac14c_1^2(X),$$ so that if $2h^1(\mathscr{O}_X)+\frac14 c_1^2(X) = 0,$ we must have that $\sigma(X) < 0.$

    Since $X$ is spin, it must be minimal. If $c_1^2(X) <0,$ then $X$ is an irrational ruled surface \cite[Theorem on p.558]{GH}, which has signature 0. So, we may assume that $c_1^2(X) = h^1(\mathscr{O}_X) = 0.$ By the classification of surfaces, we may conclude that $X$ is either $\mathbb{CP}^2,$ a Hirzebruch surface, an Enriques surface, or elliptic \cite[p. 520, p. 590]{GH}. It cannot be $\mathbb{CP}^2$ or Enriques since these are not spin. It cannot be a Hirzebruch surface since all such surfaces have signature 0. Therefore, $X$ must be elliptic as desired.
\end{proof}

Therefore, the bound in Corollary \ref{su(2)} is sharp precisely for all K3 surfaces, as we have seen, and projective elliptic surfaces with $h^1(\mathscr{O}_X)= 0$ and $K_X$ even. We arrive at the following explicit computation of $n_2(X)$:

\begin{theorem}\label{elliptic}
    If $(X,H)$ is a projective elliptic surface with $h^1(\mathscr{O}_X) = 0$ and $K_X$ even, then $$n_2(X) = n_2(X, H) = 2p_g(X) + 2.$$
\end{theorem}

\begin{proof}
    $X$ must be minimal because it is spin. The signature of a minimal elliptic surface $X$ with $h^1(\mathscr{O}_X) = 0$ is given by $$\sigma(X) = -8-8p_g(X)$$ (see for instance \cite[Theorem 24, Chapter 7]{Friedman}.) If $X$ is spin, we must have that $p_g$ is odd by Rokhlin's Theorem. Hence, combining Gieseker's Theorem with Corollary \ref{su(2)}, we arrive at the desired computation for $n_2(X, H).$
\end{proof}

Note that Qi Xia proved that if $(X, H)$ is a simply connected polarized elliptic surface without multiple fibers, then $n_2(X, H) \geq 2p_g-\epsilon,$ where $\epsilon = 1$ if $p_g$ is even and $\epsilon = 0$ if $p_g$ is odd \cite{Xia}. By Kodaira's canonical bundle formula (see for instance \cite[Theorem 15, Chapter 7]{Friedman}) a relatively minimal elliptic surface with no multiple fibers is spin if and only if $p_g$ is odd. So, in the case where $X$ is relatively minimal with $p_g$ odd, we improve Xia's bound from $n_2(X, H) \geq 2p_g$ to the sharp bound $n_2(X, H) \geq 2p_g+2.$

\begin{remark}\label{donaldsonpoly}
    Let $X$ be a simply connected elliptic surface without multiple fibers. For each positive integer $n,$ there exists such an $X$ with $\chi(\mathscr{O}_X) = n.$ As discussed above, $X$ is spin if and only if $\chi(\mathscr{O}_X)$ is even. Fintushel and Stern computed, building on work of Kronheimer and Mrowka, that the $SU(2)$ Donaldson series $D_X$ of $X$ is given by: $$D_X = \exp(Q_X/2)\sinh^{\chi(\mathscr{O}_X)-2}(f),$$ where $f\in H_2(X;\mathbb{Z})$ is the homology class of a general fiber \cite[Theorem 3.1]{FintushelStern}. Now, assume that $\chi(\mathscr{O}_X)$ is divisible by 4. Then the right-hand side is a power series that has nonzero degree $4\delta$-part for each positive integer $\delta.$ Hence, the degree $4\delta$ Donaldson polynomial is nonzero whenever $\delta \geq 1.$ This polynomial corresponds to the moduli space of anti-self-dual Yang-Mills connections with second Chern class given by $c_2 = \delta + \frac34\chi(\mathscr{O}_X).$ When $\delta<\frac54\chi(\mathscr{O}_X),$ Corollary \ref{su(2)} implies this moduli space consists only of obstructed connections. So, we have furnished infinitely many diffeomorphism types of algebraic surfaces with nonvanishing Donaldson polynomials corresponding to moduli spaces which are obstructed for all K\"ahler metrics.
\end{remark}

\subsection{Lower Bounds}

We now compare our lower bound on $n_2(X)$ with the following lower bound on $n_2(X, \omega)$ coming from a more elementary construction, but requiring stability. Our lower bound is particularly interesting when it represents a stronger bound than the below bound, and we will find that this is often the case.

\begin{theorem}
    Let $(X, \omega)$ be a K\"ahler surface. Then, $$n_r(X,\omega)\geq \frac{r^2 - 1}{2r}\chi(\mathscr{O}_X).$$
\end{theorem}

There are two ways to prove this; we demonstrate both here.

\begin{proof}[Proof using simplicity]
    Assume that a $\omega$-stable, good vector bundle $\mathscr{E}\to X$ exists with rank $r, c_1(\mathscr{E}) = 0.$ We will bound $c_2(\mathscr{E}).$ We can use the fact that the bundle $\mathscr{E}$ is \emph{simple}, i.e. that $H^0(\End_0\mathscr{E}) = 0.$ Since $$H^2(\End_0\mathscr{E}) = H^0(\End_0\mathscr{E}\otimes K)=0$$ by Serre duality as well, we must have that $\chi(\End_0\mathscr{E}) = -h^1(\End_0\mathscr{E}).$ Also, the additivity and the multiplicativity of the Chern character allow us to compute that: $$\text{ch}(\End_0\mathscr{E}) = (r^2 -1)-2rc_2(\mathscr{E}).$$ Hence, the Hirzebruch-Riemann-Roch Theorem gives: $$(r^2 - 1)\chi(\mathscr{O}_X)-2rc_2(\mathscr{E}) = \chi(\End_0\mathscr{E}) = -h^1(\End_0\mathscr{E}) \leq 0,$$ which tells us: $$c_2(\mathscr{E}) \geq \frac{r^2-1}{2r}\chi(\mathscr{O}_X),$$ as claimed.
\end{proof}

\begin{proof}[Proof using the virtual dimension]
   Assume that an $\omega$-stable, good bundle $\mathscr{E}$ exists. Then, the virtual dimension of the moduli space of $\omega$-stable rank $r$ bundles with $c_1 = 0, c_2 = c_2(\mathscr{E})$ is nonnegative. As calculated by Atiyah, Hitchin, and Singer (see \cite[Table 8.1]{Atiyah-Hitchin-Singer} or also \cite[p. 137]{DK},) this complex dimension is given by: $$2rc_2(\mathscr{E}) - (r^2 - 1)\chi(\mathscr{O}_X).$$ Hence, the nonnegativity of this virtual dimension amounts to the inequality: $$c_2(\mathscr{E})\geq \frac{r^2 - 1}{2r}\chi(\mathscr{O}_X),$$ as desired.
\end{proof}

Hence, the lower bound from Corollary \ref{realstable} represents a non-trivial improvement upon this classical bound on $n_2(X, \omega)$ whenever $$c_1^2(X) < 4(1+\frac1{r^2})\chi(\mathscr{O}_X).$$ In particular, our lower bound on $n_2(X)$ is stronger than this bound whenever $$c_1^2(X) < 5\chi(\mathscr{O}_X).$$

\begin{remark}
    Our result only applies to minimal complex surfaces since non-minimal complex surfaces are never spin. If $c_1^2(X) = 0$ then, as we have already seen, our lower bound on $n_2(X)$ is often sharp. For $X$ a minimal surface of general type, with $c_1^2(X) > 0$, the Noether and Bogomolov-Miyaoka-Yau inequalities constrain $\chi(\mathscr{O}_X)$ and $c_1^2(X)$ as follows: $$2\chi(\mathscr{O}_X)-6\leq c_1^2(X)\leq 9\chi(\mathscr{O}_X),$$ and in fact there are several examples of minimal surfaces of general type (to be seen below!) with $c_1^2(X) < 5\chi(\mathscr{O}_X)$.
\end{remark}

\subsubsection{Branched Covers of $\mathbb{CP}^2$}

The lower bound in Corollary \ref{su(2)} is particularly interesting for branched covers of $\mathbb{CP}^2.$ We will illustrate this numerically. Fix an integer $p$ and let $X$ be a two-fold cover of $\mathbb{CP}^2$ branched over a curve of degree $2p$. The surface $X$ is spin whenever $p$ is odd, and satisfies: $$\chi(\mathscr{O}_X) = \frac12p^2 - \frac32 p +2, ~ c_1^2(X) = 2p^2 - 12p + 18$$ (see for instance \cite[pp. 14-15]{DK}.) Hence, in this case, $$n_2(X) \geq \frac12 p^2  - \frac12.$$ The classical bound from simplicity is that $$n_2(X,\omega) \geq \frac38p^2 - \frac98 p +\frac32,$$ for $\omega$ a K\"ahler form on $X.$ Combining our bound with Gieseker's upper bound, we see that when $X$ is algebraic with a polarization $H$: $$ \frac12 p^2  - \frac12\leq n_2(X)\leq n_2(X, H) \leq p^2 - 3p +4.$$ When $p=3,$ $X$ is a K3 surface, and we see once again that $n_2(X)=n_2(X, H) =4.$ When $p=5,$ $X$ is a spin Horikawa surface \cite{Horikawa}, and we see that $12\leq n_2(X)\leq n_2(X, H)\leq 14.$

\subsubsection{Hypersurfaces in $\mathbb{CP}^3$}

The lower bound on $n_2(X)$ in Corollary \ref{su(2)} is also an improvement from the classical bound when $X$ is a degree $d$ hypersurface in $\mathbb{CP}^3,$ for $d\in\{4, 6,8,10, 12, 14\},$ which is spin since it has even degree. In this case $$\chi(\mathscr{O}_X) = \frac16 d^3-d^2+\frac{11}6d, ~ \sigma(X) = -\frac13d^3 + \frac43 d$$ (see for instance \cite[pp. 12-13]{DK}.) For $d = 4,$ $X$ is a K3 surface and our bound is sharp but already known. For $d>4,$ our bound appears to be new, and we make the following table illustrating the classical lower bound assuming stability, the lower bound from Corollary \ref{su(2)}, and the upper bound from Gieseker assuming algebraicity (where we round up the classical bound as necessary): 
\medskip

\begin{tabular}{ c | c | c }
Degree & Classical lower bound & Corollary \ref{su(2)} and Gieseker bound \\
\hline
6 & 9 & $16\leq n_2 \leq 22$ \\
\hline
8 & 27 & $40\leq n_2 \leq 72$ \\
\hline
10 & 64 & $80\leq n_2 \leq 170$ \\
\hline
12 & 125 & $140\leq n_2 \leq 332$ \\
\hline
14 & 216 & $224\leq n_2 \leq 574$ \\
\end{tabular}

\section{Applications to Geometric Analysis}

Taubes characterized rescaled limits of divergent sequences of solutions to the multi-spinor Seiberg-Witten equations; such rescaled limits are known as singular Fueter sections \cite{Taubes16}. In this final section of our paper, we also use the finite-dimensional approximation method on the multi-spinor Seiberg-Witten equations to prove a topological condition guaranteeing the existence of singular Fueter sections.

\subsection{Singular Fueter sections and Taubes' Compactness Theorem}

We first review Taubes' Compactness Theorem for solutions to the multi-spinor Seiberg-Witten equations. For the reader's convenience, we repeat the review in \cite[Section 4.2]{Thakar26} here. We will begin with the data given in the geometric set-up of Section 2.1 of this paper, as well as the notation $R_{spin}$ to denote the curvature of the metric connection on the spinor bundle of $X.$

\begin{definition}\label{limiting-config}
    A \emph{singular Fueter section} $(A, \Phi, Z)$ on $X$ for the metric $g$ and auxiliary connection $B$ consists of the following data:\begin{enumerate}
        \item $Z\subset X$ a closed, nowhere dense subset.
        \item A unitary connection $A$ on the line bundle $S|_{X - Z}$.
        \item A section $\Phi\in \Gamma_{X-Z}(W^+\otimes E)$ satisfying: \begin{enumerate}
            \item $||\Phi||_{L^2(X-Z)} = 1.$
            \item The norm $|\Phi|$ extends to a H\"older continuous function on $X$ which is zero on $Z.$
            \item $|\nabla_{AB}\Phi|\in L^2(X-Z).$
            \item $D_{A\otimes B}^+\Phi = 0$
            \item $\mu(\Phi, \Phi) = 0$
            \item $F_{A} = - |\Phi|^{-2}\langle\Phi, (R_{spin}+F_B)\Phi\rangle$
        \end{enumerate}
    \end{enumerate}
\end{definition}

\begin{definition} \label{convergence}
    A sequence $(A_i,\Phi_i)$ is said to \emph{converge} to a singular Fueter section $(A, \Phi, Z)$ if: \begin{enumerate}
        \item On compact subsets of $X-Z,$ \begin{enumerate}
            \item $A_i\to A$ in the weak $L^2_1$ topology, and
            \item $||\Phi_i||_{L^2}^{-1}\Phi_i\to \Phi$ in the $L^2_2$-topology.
        \end{enumerate}
        \item On $X,$ \begin{enumerate}
            \item $||\Phi_i||_{L^2}^{-1}|\Phi_i|\to |\Phi|$ in $L^2_1$ and $C^0$, and
            \item $||\Phi_i||_{L^2}^{-1}|\nabla_{A_iB}\Phi_i|\to |\nabla_{AB}\Phi|$ in $L^2$.
        \end{enumerate}
    \end{enumerate}
\end{definition}

The following beautiful theorem of Taubes exhibits an alternative: either sequences of solutions to the multi-spinor Seiberg-Witten equatiosn converge to another solution, or they converge in the way defined above to a singular Fueter section:

\begin{theorem}[{\cite{Taubes16}}]\label{taubes}
    Let $(A_i,\Phi_i)$ be a sequence of solutions to the equation $\mathcal{S}(A_i, \Phi_i) = 0$ with respect to the metric $g$ and auxiliary connection $B.$ Then, \begin{enumerate}
        \item If $\limsup\limits_{i\to\infty}||\Phi_i||_{L^2}<\infty,$ there is a subsequence of $(A_i,\Phi_i)$ which converges modulo gauge to a solution $(A, \Phi)$ of $\mathcal{S}(A, \Phi) = 0$ with respect to $(g, B).$
        \item If $\limsup\limits_{i\to\infty}||\Phi_i||_{L^2}=\infty,$ there is a subsequence of $(A_i,\Phi_i)$ which converges modulo gauge to a singular Fueter section $(A, \Phi, Z).$
    \end{enumerate}
\end{theorem}

\subsection{Existence Theorem}

We will now show how Theorem \ref{Furuta} can be used to show topological conditions guaranteeing the existence of singular Fueter sections and solutions to the multi-spinor Seiberg-Witten equations. The power of the below theorem is that once the topological criterion is satisfied, the existence will be guaranteed for every Riemannian metric. Hene, we can view it as a 4-dimensional analog of existence theorems for $\mathbb{Z}/2$-harmonic 1-forms on 3-manifolds that hold true for every Riemannian metric on a fixed smooth manifold, coming from an interpretation as limits of unbounded sequences of solutions to a gauge-theoretic equation that has non-compact moduli in all metrics. This is best exemplified by Taubes' seminal theorem that $\mathbb{Z}/2$-harmonic 1-forms exist on every Riemannian 3-manifold with non-compact $SL(2;\mathbb{C})$ character variety \cite{TaubesSL2C}.

If $V$ is an $SO(r)$-bundle, and $B$ is a $SO(r)$-connection, then we can view $B$ as a connection in the $SU(r)$-bundle $V\otimes \mathbb{C}$ via the natural inclusion $SO(r)\subset SU(r)$. 

\begin{theorem}\label{existence}
    Suppose $X$ is a spin, closed, smooth, oriented 4-manifold and $V\to X$ is a $SO(r)$-bundle with complexification $E$. If \begin{equation*}\label{eq: singularfuetercondition} p_1(V) > \frac{r\sigma(X)}{8}+\max\left(0, b_+(X) - 1\right),\end{equation*} then for every Riemannian metric $g$ and $SO(r)$-connection $B$ on $V,$ the moduli space $\mathcal{M}(g, B)$ is non-empty, it contains an irreducible solution, and there exists a singular Fueter section on $X$ with respect to $g$ and $B.$
\end{theorem}

\begin{proof}
    Fix a metric $g$ and $SO(r)$-connection $B$ on $V$, which we can view as a connection on $E$ by the remark above. It has an obvious real structure on it for which $B$ is compatible. Therefore, if the moduli space $\mathcal{M}(g, B)$ is compact, Theorem \ref{Furuta} gives that $$ -c_2(E) \leq\frac{r\sigma(X)}{8}+\max\left(0, b_+(X)- 1\right).$$ Since $c_2(E) = -p_1(V),$ the bound from Theorem \ref{Furuta} contradicts our assumption. Therefore, $\mathcal{M}(g, B)$ is not compact, so in particular it is non-empty. Theorem \ref{taubes} now tells us that there exists a singular Fueter section with respect to $g$ and $B.$ The moduli space must contain an irreducible solution since the space of reducible solutions is compact (see the proof of Theorem \ref{compactness}.)
\end{proof}

\begin{remark}
    This theorem is particularly interesting in the case when the signature $\sigma(X)$ is negative. For, even if $X$ has zero Seiberg-Witten invariant, then this argument shows $\mathcal{M}(g, B)$ is nonetheless non-empty once $r$ is sufficiently large relative to $-p_1(V)$. For example, $X = K3\# K3$ has vanishing Seiberg-Witten invariants since it is a connected sum of two 4-manifolds with $b_+\geq1,$ however it has irreducible solutions to the two-spinor Seiberg-Witten equations for every Riemannian metric whenever $-p_1(V)\leq 3$.
\end{remark}

\appendix
\section{Bryan's proof of Proposition \ref{ktheory}}

We present Bryan's beautiful proof of the following proposition, originally due to Furuta:

\begin{proposition}
    Suppose there exist nonnegative $m, n,$ positive integers $\beta, \kappa,$ a real number $R>0,$ and a $\Pin(2)$-equivariant map $$f:\mathbb{H}^{\kappa+m}\oplus \widetilde{\mathbb{R}}^n\to \mathbb{H}^{m}\oplus \widetilde{\mathbb{R}}^{\beta+n}$$ that has no zeroes on $S_R(\mathbb{H}^{\kappa+m}\oplus \widetilde{\mathbb{R}}^n).$ Then, $$2\kappa+1\leq \beta.$$
\end{proposition}

\subsection{K-theory and the tom Dieck character formula}
Let $\Gamma$ be a compact Lie group with representation ring $R(\Gamma)$, and let $V$ and $W$ be complex representations. Then, we define the equivariant $K$-theory $K_\Gamma(V)$ to be the relative equivariant $K$-theory $K_\Gamma(BV, SV)$ where $BV$ is a $\Gamma$-invariant closed neighborhood of the origin in $V$ and $SV$ is its boundary. 

The equivariant Thom isomorphism in $K$-theory shows that $K_\Gamma(V)$ is a free $R(\Gamma)$-module generated by the Bott class $\Lambda(V)\in K_\Gamma(V)$ \cite[Theorem 4.3]{equivariantthom}. A $\Gamma$-equivariant function $F:(BV, SV)\to (BW, SW)$ gives us a map on $K$-theory $F^*:K_\Gamma(W)\to K_\Gamma(V).$ Therefore, to a $\Gamma$-equivariant function $F:(BV, SV)\to (BW, SW)$ we may associate the \emph{$K$-theoretic degree} $\alpha_F \in R(\Gamma)$, which is the unique element of the representation ring satisfying $$\alpha_F\Lambda(V) = F^*(\Lambda(W)).$$ The tom Dieck character formula provides a formula to compute $\alpha_F.$ To introduce this, we will first discuss some notation. For each $\gamma\in \Gamma,$ let $V^\gamma, W^\gamma$ denote the fixed subspaces of $V, W$ by $\gamma.$ Since $F$ is equivariant, the restriction $F|V^\gamma:V^\gamma\to W^\gamma$ is well-defined and we denote by $\deg(F, \gamma)$ the ordinary topological degree of the map $$F:(BV^\gamma, SV^\gamma)\to (BW^\gamma, SW^\gamma).$$ Following \cite{tomDieck}, for each $x\in R(\Gamma)$ we define $\lambda_{-1}(x)$ to be the alternating sum of exterior powers $\lambda_{-1}(x) = \sum_k (-1)^k\Lambda^kx$. We have:

\begin{theorem}\cite[Proposition 9.7.2]{tomDieck}
    For each $\gamma \in \Gamma,$ we have the following equality of traces: $$\tr_\gamma(\alpha_F) = \deg(F, \gamma)\tr_\gamma\left(\lambda_{-1}\left((W^\gamma)^\perp - (V^\gamma)^\perp\right)\right).$$
\end{theorem}

\subsection{Irreducible representations of $\Pin(2)$}
First, we note that the complex irreducible representations of $S^1$ are the characters of $S^1,$ and we denote the character $z\mapsto z^k$ by $\theta_k.$

The complex irreducible representations $\Pin(2)$ are the trivial representation $\mathbb{C},$ the representation $\widetilde{\mathbb{C}} := \widetilde{\mathbb{R}}\otimes_\mathbb{R}\mathbb{C},$ and 2-dimensional representations $\mathbb{H}_k,$ for $k\in\{1,2,\dots\}.$ The restriction of each $\mathbb{H}_k$ to $S^1\subset \Pin(2)$ is given by $\theta_k+\theta_{-k}.$ The representation $\mathbb{H}_1$ is the quaternion representation $\mathbb{H}_1 = \mathbb{H}.$

\subsection{The proof}
\begin{proof}[Proof of Proposition \ref{ktheory}]
    In this section, we will use $V_\mathbb{R}$ and $W_\mathbb{R}$ to denote the domain and codomain, respectively, of $f$ as real $\Pin(2)$-modules. We define the complexifications $V = V_\mathbb{R}\otimes_\mathbb{R}\mathbb{C}$ and $W= W_\mathbb{R}\otimes_\mathbb{R}\mathbb{C}$. These are complex $\Pin(2)$-modules, with irreducible decomposition $$V = \mathbb{H}^{2\kappa+2m}\oplus \widetilde{\mathbb{C}}^n,~ W=\mathbb{H}^{2m}\oplus \widetilde{\mathbb{C}}^{\beta+n}.$$ We may define a function $\tilde{f}:V\to W$ as: $$\tilde{f}(u\otimes 1+v\otimes i) = f(u)\otimes 1+f(v)\otimes i.$$ Define $BV$ to be the $\Pin(2)$-invariant polyball $$BV := \{u\otimes 1+v\otimes i\in V | ~||u||_V\leq R, ||v||_V\leq R\},$$ and denote by $SV$ its boundary. Denote by $BW$ the unit ball in $W$ and $SW$ its boundary. Clearly $\tilde{f}|SV$ has no zeroes, so the map $$\frac{\tilde{f}}{||\tilde{f}||_W}:SV\to SW$$ is well-defined and $\Pin(2)$-equivariant, and we may define a $\Pin(2)$-equivariant map $F:(BV, SV)\to (BW, SW)$ by coning off the map $\frac{\tilde{f}}{||\tilde{f}||_W}.$ 

    We know that the $K$-theoretic degree can be written as $$\alpha_F = a[\mathbb{C}]+b[\widetilde{\mathbb{C}}]+\sum_{k=1}^\infty c_k[\mathbb{H}_k],$$ for integers $a, b, c_k.$ Let $\varphi\in S^1$ be any topological generator. Since $\varphi$ acts trivially on $\widetilde{\mathbb{C}}$ but non-trivially on $\mathbb{H},$ we have that $\dim V^\varphi - \dim W^\varphi = \beta>0,$ hence $\deg(F, \varphi) = 0.$ So, from the tom Dieck character formula, we conclude that $$0 = \tr_\varphi(\alpha_F) = a+b+\sum_{k=1}^\infty c_k(\varphi^k+\varphi^{-k}).$$ Since this is true for any topological generator $\varphi\in S^1,$ we conclude $a=-b, c_k = 0$ for all $k.$ 

    Now, we compute the trace of $j$ on $\alpha_F.$ Since $j$ acts non-trivially on both $\widetilde{\mathbb{C}}$ and $\mathbb{H},$ we have that $\dim V^j = \dim W^j = 0,$ so $\deg(F, j) = 1.$ Therefore, using the facts that $\lambda_{-1}\left(a_1r_1+a_2r_2\right) = (\lambda_{-1}r_1)^{a_1} (\lambda_{-1}r_2)^{a_2}$ and $\tr_J\mathbb{H} = 0,$ we have: \begin{equation*}
        \begin{split}
            2a &= a + (-b) \\
            &= \tr_j(\alpha_F) \\
            &= \deg(F, j)\tr_j\left(\lambda_{-1}\left((W^j)^\perp - (V^j)^\perp\right)\right) \\
            &= \tr_j\left(\lambda_{-1}\left(W - V\right)\right) \\
            &= \tr_j\left(\lambda_{-1}\left(\beta[\widetilde{\mathbb{C}}] - 2\kappa[\mathbb{H}]\right)\right) \\
            &= \tr_j\left((\lambda_{-1}[\widetilde{\mathbb{C}}])^\beta(\lambda_{-1}[\mathbb{H}])^{-2\kappa}\right) \\
            &= \tr_j\left((\lambda_{-1}[\widetilde{\mathbb{C}}])^\beta(\lambda_{-1}[\mathbb{H}])^{-2\kappa}\right) \\
            &= \tr_j\left(([\mathbb{C}] - [\widetilde{\mathbb{C}}])^\beta(2[\mathbb{C}]-[\mathbb{H}])^{-2\kappa}\right) \\
            &= 2^{\beta - 2\kappa} \\
        \end{split}
    \end{equation*}

    Therefore $a = 2^{\beta-2\kappa-1}.$ Since $a$ is an integer, we must have $2\kappa+1<\beta.$
\end{proof}

\bibliographystyle{alpha}
\bibliography{refs2}

\end{document}